\documentstyle[12pt]{article}

\input amssym.def
\input amssym.tex

\def\ie{{\it i.e. }}
\def\eg{{\it e.g. }}
\def\={\ = \ }
\def\la{\langle}
\def\ra{\rangle}

\def\be{\setcounter{equation}{\value{theorem}} \begin{equation}}
\def\ee{\end{equation} \addtocounter{theorem}{1}}
\def\bp{{\sc Proof: }}
\def\ep{{}{\hfill $\Box$} \vskip 5pt \par}
\def\bl{\begin{lemma}}
\def\el{\end{lemma}}
\def\bt{\begin{theorem}}
\def\et{\end{theorem}}
\def\bprop{\begin{prop}}
\def\eprop{\end{prop}}

\def\H{{\cal H}}

\def\M{{\cal M}}

\def\D{{\Bbb D}}
\def\C{{\Bbb C}}
\def\B{{\Bbb B}}

\def\T{T_\Phi}

\def\la{\langle}
\def\={\ = \ }

\def\l{\lambda}
\def\L{\Lambda}

\def\rxl{R_{x,\lambda}}
\def\rxL{R_{x,\Lambda}}
\def\trxL{R_{x,\Lambda}^\sim}

\advance\topmargin by -.6in
\advance\textheight by 1.2in
\advance\oddsidemargin by -.5in
\advance\textwidth by 1.05in

\newtheorem{theorem}{Theorem} [section]
\newtheorem{prop}[theorem]{Proposition}
\newtheorem{lemma}[theorem]{Lemma}
\newtheorem{cor}[theorem]{Corollary}

\begin{document}
\setlength{\baselineskip}{21pt}
\title{ Complete Nevanlinna-Pick Kernels}
\author{Jim Agler
\thanks{Partially supported by the National Science Foundation}
\\{John E. M\raise.5ex\hbox{c}Carthy
\thanks{Partially supported by National Science Foundation
grant DMS 9531967. }}
\date{}
\\{\rm University of California at San Diego, La Jolla California 92093}
\\{\rm Washington University, St. Louis, Missouri 63130}}

\bibliographystyle{plain}
\maketitle

\begin{abstract}
We give a new treatment of Quiggin's and 
M\raise.45ex\hbox{c}Cullough's characterization of
complete Nevanlinna-Pick kernels. 
We show that a kernel has the matrix-valued
Nevanlinna-Pick property if and only if it has the vector-valued
Nevanlinna-Pick property.
We give a representation of all complete Nevanlinna-Pick kernels, and show
that they are all restrictions of a universal complete Nevanlinna-Pick
kernel.
\end{abstract}

\baselineskip = 18pt

\setcounter{section}{-1}
\section{Introduction}
Let $X$ be an infinite set, and $k$ a positive definite kernel function on
$X$, \ie for any finite collection  $x_1,\dots, x_n$ of distinct points in 
$X$, and any complex
numbers $\{ a_i \}_{i=1}^n$, the sum 
\be
\sum_{i,j=1}^n a_i \bar a_j k(x_i, x_j) \geq 0,
\label{pd}
\ee
with strict inequality unless all the $a_i$'s are $0$.
For each element $x$ of $X$, define the function $k_x$ on $X$
by $k_x(y) := k(x,y)$.
Define an inner product on the span of these functions by
$$
\la \sum a_i k_{x_i} \, , \, \sum b_j k_{y_j} \ra \=
\sum a_i \bar b_j k(x_i,y_j),
$$
and let $\H = \H_k$ be the Hilbert space obtained by completing the space of
finite linear combinations of $k_{x_i}$'s with respect to this inner
product. 
The elements of $\H$ can be thought of as functions on $X$, with the value of
$f$ at $x$ given by $\la f,k_x \ra$.

A multiplier of $\H$ is a function $\phi$ on $X$ with the property that if $f$
is in $\H$, so is $\phi f$. The Nevanlinna-Pick problem is to determine, given
a finite set $x_1, \dots, x_{n}$ in $X$, and numbers $\lambda_1,
\dots, \lambda_{n}$, whether there exists a multiplier $\phi$ of norm at
most one that interpolates the data, \ie satisfies $\phi(x_i) = \lambda_i$
for $i = 1, \dots, n$. 

If $\phi$ is a multiplier of $\H$, we shall let $T_\phi$ denote the operator
of multiplication by $\phi$. Note that the adjoint of $T_\phi$ satisfies
$T_{\phi}^\ast k_x = \overline{\phi(x)} k_x$. So if $\phi$ interpolates the data
$(x_i,\lambda_i)$, then the $n$-dimensional space spanned by $\{ k_{x_i} :
1 \leq i \leq n \}$ is left invariant by $T_{\phi}^\ast$, and on this subspace
the operator $T_{\phi}^\ast$ is the diagonal
\be
\label{di}
\left(\matrix{
\overline{\l_1} & & \cr
& \ddots & \cr
&& \overline{\l_{n}} \cr}
\right)
\ee
with respect to the (not necessarily orthonormal) basis $\{ k_{x_i} \}$.
For a given set of $n$ data points $(x_1,\lambda_1), \dots, (x_n,\l_n)$,
let $\rxl$ be the operator in (\ref{di}), \ie the operator that sends
$k_{x_i}$ to $ \overline{\l_i} k_{x_i}$.
A necessary condition to solve the Nevanlinna-Pick problem is that
the norm of $\rxl$ be at most $1$; the kernel $k$ is called a
{\it Nevanlinna-Pick kernel} if this necessary condition is also always 
sufficient. 

Notice that $\rxl$ is a contraction on $sp\{ k_{x_i} :
1 \leq i \leq n \}$ if and only if $(1 - \rxl^\ast \rxl)$ is positive on
that space. As
$$
\la (1 - \rxl^\ast \rxl) \sum_{i=1}^n a_i k_{x_i} ,
\sum_{j=1}^n a_j k_{x_j} \ra \=
\sum_{i,j=1}^n a_i \overline{a_j} ( 1 - \l_j \overline{\l_i})
\la k_{x_i} , k_{x_j} \ra ,
$$
it follows that the contractivity of $\rxl$ on $sp\{ k_{x_i} :
1 \leq i \leq n \}$ is equivalent to the positivity of the $n$-by-$n$ matrix
\be
\label{eq:pos}
\left( ( 1 - \l_j \overline{\l_i})
\la k_{x_i} , k_{x_j} \ra \right)_{i,j=1}^n .
\ee

The classical Nevanlinna-Pick theorem asserts that the 
Szeg\"o kernel 
$$
k(x,y) = \frac{1}{1- \bar x y}
$$
on the unit disk is a Nevanlinna-Pick kernel. The condition is normally
stated in terms of  the positivity of 
(\ref{eq:pos}), but as we see that is equivalent to the contractivity of
(\ref{di}).

The matrix-valued Nevanlinna-Pick problem is as follows.
Fix some auxiliary Hilbert space, which for notational convenience we shall 
assume to be the finite-dimensional space $\C^\nu$. 
The tensor product $\H \otimes \C^\nu $ 
can be thought of as a space of vector valued functions on $X$.
A multiplier of $\H \otimes \C^\nu $ is now a $\nu$-by-$\nu$ matrix valued
function $\Phi$ on $X$ with the property that whenever
$$
\left(\matrix{
f_1 \cr
\vdots  \cr
f_\nu \cr}
\right)
\ \in \
\H \otimes \C^\nu
$$
then
$$
\Phi
\left(\matrix{
f_1 \cr
\vdots  \cr
f_\nu \cr}
\right)
\ \in \
\H \otimes \C^\nu
.
$$
The matrix Nevanlinna-Pick problem is to determine,
given points $x_1,\dots,x_{n}$ and matrices 
$\L_1,\dots,\L_{n}$, whether there is a multiplier $\Phi$ of norm at most
one that interpolates: $\Phi(x_i) = \Lambda_i$.

Fix a (not necessarily orthonormal)
basis $\{ e^\alpha \}_{\alpha=1}^\nu$ for $\C^\nu$.
As before, 
$$
\T^\ast k_{x} \otimes v  \= k_x \otimes \Phi(x)^\ast v ,
$$
so if $\M$ is the span of $\{ k_{x_i} \otimes e^\alpha : 1 \leq i \leq n,\, 
1\leq \alpha
\leq \nu \}$, a necessary condition for the Nevanlinna-Pick problem to have a
solution is that the $n\nu$-by-$n\nu$ matrix
\be
\label{17}
\rxL :\ k_{x_i} \otimes e^\alpha \mapsto k_{x_i} \otimes \L_i^\ast e^\alpha 
\ee
be a contraction. We shall call the kernel $k$ a {\it complete Nevanlinna-Pick
kernel} if, for all finite $\nu$ and all positive $n$,
the contractivity of $\rxL$ is also a sufficient condition to extend
$\Phi$ to a multiplier of all of $\H \times \C^\nu$ of norm at most one.
 
In Section~\ref{sec1} we
give a classification of  all complete Nevanlinna-Pick kernels.
This was originally done by S.~M\raise.45ex\hbox{c}Cullough
in \cite{mccul92} in the context of the Carath\'eodory interpolation problem.
The Nevanlinna-Pick problem was studied by P.~Quiggin, 
who in \cite{qui93}
established the sufficiency
of the condition in Theorem~\ref{thm}, and in \cite{qui94}
established the necessity.

In Section~\ref{sec4} we show that if a kernel has the
Nevanlinna-Pick property for row vectors of length $\nu$, then it has the
Nevanlinna-Pick property for $\mu$-by-$\nu$ matrices for all $\mu$.
In particular, having the vector-valued Nevanlinna-Pick property 
is equivalent to having the complete Nevanlinna-Pick property.

In Section~\ref{sec2}, we show that all complete Nevanlinna-Pick kernels have
the form
$$
k(x,y) \= \frac{\overline{\delta(x)} \delta(y)}{1 - F(x,y)}
$$
where $\delta$ is a nowhere vanishing function and $F:X \times X \to \D$
is a positive semi-definite function.

In Section~\ref{sec3} we introduce the universal complete Nevanlinna-Pick
kernels $a_m$ defined on the unit ball $\B_m$ of an $m$-dimensional Hilbert space
($m$ may be infinite) by
$$
a_m(x,y) \= \frac{1}{1-\la x,y \ra}.
$$
These kernels are universal in the sense that, up to renormalization, every 
complete Nevanlinna-Pick kernel is just the restriction of an $a_m$ to a
subset of $\B_m$.

\section{Characterization of Complete Nevanlinna-Pick kernels}
\label{sec1}

To simplify notation, we shall let $k_i$ denote $k_{x_i}$, and $k_{ij}$
denote $\la k_i, k_j \ra = k(x_i,x_j)$.
First we want a lemma that says that we can break $\H$ up into summands on
each of which $k$ is irreducible, \ie $k_{ij}$ is never $0$.
For convenience, we shall defer the proof of the lemma until after the proof
of the theorem.

\begin{lemma}
\label{lem1}
Suppose $k$ is a Nevanlinna-Pick kernel on the set $X$.
Then $X$ can be partitioned into disjoint subsets $X_i$ such that if two
points $x$ and $y$ are in the same set $X_i$, then $k(x,y) \neq 0$; and
if $x$ and $y$ are in different sets of the partition, then $k(x,y) = 0$.
\end{lemma}

A reducible kernel will have the (complete) Nevanlinna-Pick
property if and only if each irreducible piece does, so we shall assume $k$
is irreducible.

\bt
\label{thm}
A necessary and sufficient condition 
for an irreducible kernel
$k$ to be a complete Nevanlinna-Pick
kernel is that, for any finite set $\{ x_1, \dots, x_n \}$ of $ n$
distinct
elements of $X$, the $(n-1)$-by-$(n-1)$ matrix
\be
\label{sm}
 F_n \= \left( 1 - \frac{k_{in} k_{nj}}{k_{ij}k_{nn}} \right)_{i,j=1}^{n-1}
\ee
is positive semi-definite.
\et
\bp
Let $x_1, \dots, x_{n-1}$  and $\L_1, \dots , \L_{n-1}$ be chosen, 
let $\M$ be the span of $\{ k_i \otimes e^\alpha : 1 \leq i \leq n-1, 1 \leq \alpha 
\leq
\nu \}$, and define $\rxL$ on $M$ by (\ref{17}).
The operator  $\rxL$ is a contraction if and only if 
$I - \rxL^\ast \rxL \geq 0$. 
Calculate 
\be
\label{21}
\la  (I - \rxL^\ast \rxL ) \sum_{i,\alpha} a_i^\alpha k_i \otimes e^\alpha \, , \,
\sum_{j,\beta} a_j^\beta k_j \otimes e^\beta \ra 
\=
\sum_{i,\alpha,j,\beta} a_i^\alpha \bar a_j^\beta  k_{ij} (\la e^\alpha, e^\beta \ra -
\la \L_j \L_i^\ast e^\alpha, e^\beta \ra )
\ee

A necessary and sufficient condition 
to be able to find a matrix $\L_n$ so that the extension
$\trxL$ of $\rxL$ that sends $\ k_{x_n} \otimes e^\alpha $
to $ k_{x_n} \otimes \L_n^\ast e^\alpha$ for each $\alpha$
has the same norm as $\rxL$ 
is: 
whenever $\L_1, \dots, \L_{n-1}$ are chosen so that
\be
\label{eq:87}
I - \rxL^\ast \rxL \ \geq \ 0
\ee
on $\vee\{ k_i \otimes e^\alpha : 1 \leq i \leq n-1, 1 \leq
\alpha \leq \nu \}$, then 
\be
\label{eq:88}
P - (P \trxL P)^\ast (P \trxL P) \ \geq \ 0,
\ee
where $P$ is the orthogonal projection from 
$\vee\{ k_i \otimes e^\alpha : 1 \leq i \leq n, 1 \leq \alpha \leq \nu \}$
onto the orthogonal complement of 
$\vee\{ k_n \otimes e^\alpha : 1 \leq \alpha \leq \nu \}$.
(This was first proved in \cite{ag1} in the scalar case, and
a proof of the matrix case is given in \cite{agmc_loc}.
Notice that (\ref{eq:88}) does not depend on the choice of $\Lambda_n$.
We use $\vee$ to denote the closed linear span of a set of vectors.)

That (\ref{eq:87}) always implies (\ref{eq:88}) for any choice of $x$ and $\Lambda$
is not only necessary, but also sufficient for $k$ to be a complete
Nevanlinna-Pick kernel. 
Sufficiency is proved by an inductive argument that if one can always extend
a multiplier defined on a finite set to any other point without increasing
the norm, then one can extend the multiplier to all of $X$. In the absence of
any {\it a priori} simplifying assumptions about the multiplier algebra of
$\H$ being large, the proof of this inductive argument is subtle, and is
originally due to Quiggin 
\cite[Lemma 4.3]{qui93}.


Using the fact that 
$$
P( k_i \otimes e^\alpha) \= (k_i - \frac{k_{in}}{k_{nn}} k_n ) \otimes e^\alpha,
$$
we can calculate that
$$
\la (P - (P \trxL P)^\ast (P \trxL P)
\sum_{i,\alpha} a_i^\alpha k_i \otimes e^\alpha \, , \,
\sum_{j,\beta} a_j^\beta k_j \otimes e^\beta \ra
$$
equals
\be
\label{23}
\sum_{i,\alpha,j,\beta} a_i^\alpha \bar a_j^\beta  k_{ij}  \left(
1- \frac{k_{in}k_{nj}}{k_{ij}k_{nn}} \right) 
\left[ \la e^\alpha, e^\beta \ra - \la \L_j \L_i^\ast e^\alpha, e^\beta \ra  \right]
\ee

Comparing (\ref{21}) and (\ref{23}), we see that we want that whenever the matrix
whose $(i,\alpha)^{th}$ column and $(j,\beta)^{th}$ row is given by
\be
\label{25}
k_{ij} (\la e^\alpha, e^\beta \ra -
\la \L_j \L_i^\ast e^\alpha, e^\beta \ra )
\ee
is positive, then the Schur product of this matrix with 
$F_n \otimes J$ is positive, where $J$ is the $\nu$-by-$\nu$ matrix all of
whose entries are $1$.
As the Schur product of two positive matrices is positive, the positivity of
(\ref{sm}) is immediately seen to be a sufficient condition for $k$ to be a
complete Nevanlinna-Pick kernel.

We shall prove necessity by induction on
$n$. The case $n=2$ holds by the Cauchy-Schwarz inequality.
So assume that $F_{n-1}$ is positive, and we shall prove that $F_n$ is positive.

Note first the sort of matrices that can occur in (\ref{25}). 
For each $i$ and $\alpha$, one can choose the vector $\Lambda_i^\ast
e^\alpha$ arbitrarily. In particular, let $G$ be any positive $(n-1)$-by-$(n-1)$
matrix, let $\varepsilon > 0$, and choose
$\{e^\alpha\}$ so that $\la e^\alpha, e^\beta \ra =  \varepsilon
\delta_{\alpha, \beta} + 1$.
Choose vectors $v_i$ so that $\la v_i, v_j \ra = G_{ij}$.
Let $\nu = n-1$, 
and choose $\Lambda_i^\ast$ to be the rank one matrix that sends each
$e^\alpha$ to $v_i$. Then (\ref{25})
becomes
\be
\label{32}
k_{ij} (\varepsilon \delta_{\alpha,\beta} + 1 - G_{ij}).
\ee
We know that  $F_n$ has the property that if $G$ is a positive matrix and
the $(n-1)\nu$-by-$(n-1)\nu$ matrix
(\ref{32}) is positive, then the Schur product of $F_n\otimes J$ with
(\ref{32}) is also positive. 
Denote by $K$ the $(n-1)$-by-$(n-1)$ matrix
whose $(i,j)$ entry is $k_{ij}$, and let $\cdot$ denote Schur product.
By letting $\varepsilon$ tend to
zero, we get that whenever $G \geq 0$ and
$$
[ K \cdot ( J - G) ] \otimes J \ \geq \ 0 ,
$$
then
$$
[F_n \otimes J ] \cdot ( [ K \cdot ( J - G) ] \otimes J )  \ \geq \ 0 ,
$$
which is the same as saying
\be
K\cdot ( J - G)  \geq 0 \quad \Rightarrow \quad
F_n \cdot K\cdot ( J - G)  \geq 0 .
\label{66}
\ee

Let $L$ be the rank one positive $(n-1)$-by-$(n-1)$ matrix given by
$$
L_{ij} \= \frac{k_{i(n-1)}k_{(n-1)j}}{k_{(n-1)(n-1)}} ,
$$
and let $G$ be the matrix given by
$$
G_{ij} \= 1 - \frac{L_{ij}}{k_{ij}} .
$$
Then $G$ is the matrix that agrees with $F_{n-1}$ in the first $(n-2)$ rows and
columns, and all the entries in the $(n-1)^{st}$ row and column are zero.
Therefore $G $ is positive by the inductive hypothesis.
Moreover, $K \cdot (J-G) = L $ and so is positive.
Therefore $F_n \cdot L$ is positive. But $L$ is rank one, so 
$1/L$ (the matrix of reciprocals) is also positive, and therefore
$$
F_n \cdot L \cdot 1/L \= F_n \geq 0,
$$
as desired.
\ep

{\sc Proof of Lemma \ref{lem1}:}
Let $X_x = \{ y : k(x,y) \neq 0 \}$. We need to show that for any two points
$x$ and $y$, the sets $X_x $ and $X_y$ are either equal or disjoint.
This is equivalent to proving that if $k(x,z) \neq 0$ and $k(y,z) \neq 0$,
then $k(x,y) \neq 0$.

Assume this fails. Consider the $2$-by-$2$ matrix  $T^\ast$
defined on the linear span of $k_x$ and $k_y$ by
\begin{eqnarray*}
T^\ast k_x &\=& k_x \\
T^\ast k_y &\=& - k_y
\end{eqnarray*}
This has norm one, because $k(x,y) = 0$.
By the hypothesis that $k$ is a Nevanlinna-Pick kernel, $T^\ast$ can be
extended to the space spanned by $k_x,k_y$ and $k_z$ so that the new
operator has the same norm and has
$k_z$ as an eigenvector. But for this to hold, from equation
(\ref{23}) we would need
\be
\label{27}
\left(\matrix{
0 & 2 \cr
2 & 0 \cr}
\right)
\cdot
\left(\matrix{
k_{xx} - \frac{|k_{xz}|^2}{k_{zz}} & k_{xy} - \frac{k_{xz}k_{zy}}{k_{zz}} \cr
k_{yx} - \frac{k_{yz} k_{zx}}{k_{zz}} & k_{yy} - \frac{|k_{yz}|^2
}{k_{zz}}\cr}
\right)
\
\geq \ 0 .
\ee
But the Schur product of the
 two matrices in (\ref{27}) is zero on the
diagonal, non-zero off the diagonal, and therefore cannot be positive.
\ep

By the same argument as in the theorem, an irreducible kernel will
have the (scalar) Nevanlinna-Pick property
if and only if whenever $G$ is positive and rank one,
(\ref{66}) holds. We do not know how to classify such kernels in the sense of
Theorem \ref{thm}.

The positivity of $F_n$ can be expressed in other ways. The proof that 
$F_n$ being positive is equivalent to
$1/K$ having only one positive eigenvalue below is
due to Quiggin \cite{qui93}. 
\begin{cor}
\label{cor:3}
A necessary and sufficient condition for the irreducible kernel $k$ 
to have the complete Nevanlinna-Pick 
property is that for any finite set $x_1, \dots, x_n$, the matrix
$$
H_n \ := \
\left( \frac{1}{k_{ij}} \right)_{i,j=1}^n
$$
has exactly one positive eigenvalue (counting multiplicity).
\end{cor}
\bp
As all the diagonal entries of $H_n$ are positive, $ H_n$ must have at least
one positive eigenvalue.

The condition that $F_{n+1}$ be positive is equivalent to saying
\be
\label{38}
M_n \ := \ 
\left(
\frac{k_{n+1,n+1}}{k_{i,n+1}k_{n+1,j}} - \frac{1}{k_{ij}}
\right)_{i,j=1}^n \ \geq \ 0,
\ee
because $k_{i,n+1}k_{n+1,j}$ is rank one so its reciprocal is positive.
But (\ref{38}) says that $H_n$ is less than or equal to 
a rank one positive operator, so
has at most one positive eigenvalue.

Conversely, any symmetric matrix
$$
\left(\matrix{
A & B \cr
B^\ast & C \cr}
\right)
$$
with $C$ invertible is congruent to 
$$
\left(\matrix{
A-BC^{-1}B^\ast & 0 \cr
0 & C \cr}
\right).
$$
(The top left entry is called the Schur complement of $C$.)
Applying this to $H_n$ with $C$ the $(n,n)$ entry, we get that $H_n$ is
congruent to
$$
\left(\matrix{
-M_{n-1} & 0 \cr
0 & \frac{1}{k_{nn}} \cr}
\right).
$$
So if $H_n$ has only one positive eigenvalue,
$-M_{n-1}$ must be negative semi-definite, and therefore $F_n$ must be positive
semi-definite.
\ep

As an application of the Corollary, consider the Dirichlet space of
holomorphic functions on the unit disk with reproducing kernel
$k(w,z) = \frac{1}{\bar w z} \log \frac{1}{1 - \bar w z}$.
It is shown in \cite{ag1} that this is a Nevanlinna-Pick kernel, and in the
course of the proof it is established that $1 - 1/k$ is positive
semi-definite (because all the coefficients in the power series are
positive).
It then follows at once from Corollary \ref{cor:3} that the Dirichlet kernel is
actually a complete Nevanlinna-Pick kernel.

\section{Vector-valued Nevanlinna-Pick kernels}
\label{sec4}

Let $\M_{\mu,\nu}$ denote the $\mu$-by-$\nu$ matrices. Let us say that a
kernel $k$ has the $n$-point $\M_{\mu,\nu}$ Nevanlinna-Pick property if,
for any points $x_1, \dots, x_n$, and any matrices $\L_1, \dots, \L_n$ in
$\M_{\mu,\nu}$, there exists a multiplier $\Psi$,
$$
\Psi : \H \otimes \C^\nu \rightarrow \H \otimes \C^\mu ,
$$
such that $\Psi(x_i) = \L_i, \ 1 \leq i \leq n$, and
$$
\|T_\Psi \| = \| T_{\Psi}^* \| = \| T_{\Psi}^* |_{sp \{ k_{x_i} \otimes
\C^\mu : \ 1 \leq i \leq n \} } \| .
$$
We shall say that $k$ is a vector-valued Nevanlinna-Pick kernel if $k$ has
the $n$ point $\M_{1,\nu}$ Nevanlinna-Pick property for all $n$ and $\nu$.
\begin{theorem}
Let $\nu \geq n-1$. Then $k$ has the $n$-point $\M_{\mu,\nu}$ Nevanlinna-Pick
property for some positive integer $\mu$ if and only if it has the property
for all positive integers $\mu$.
\end{theorem}
\bp
It is clear that the $n$-point $\M_{\mu,\nu}$ Nevanlinna-Pick
property implies the $n$-point $\M_{\pi,\nu}$ Nevanlinna-Pick
property for all $\pi$ smaller than $\mu$.
So it is sufficient to prove that the $n$-point $\M_{1,\nu}$
Nevanlinna-Pick property
implies the $n$-point $\M_{\mu,\nu}$ Nevanlinna-Pick
property for all $\mu$.

As in the proof of Theorem~\ref{thm}, the kernel $k$ has the
$n$-point $\M_{\mu,\nu}$ Nevanlinna-Pick
property if and only if the positivity of the matrix
\be
\label{eq:24}
\left[ k_{ij} (\la e^\alpha, e^\beta \ra_{\C^\mu}  -
\la \L_j \L_i^\ast e^\alpha, e^\beta \ra_{\C^\nu} )
\right]_{i,j=1; \alpha,\beta = 1}^{i,j=n; \alpha,\beta = \mu}
\ee
implies the positivity of the Schur product of 
(\ref{eq:24}) with $F_{n+1} \otimes J_\mu$.
Again, as in the proof of Theorem~\ref{thm}, this implies that whenever
$K \cdot (J_n - G)$ is positive, 
then so is
$F_{n+1} \cdot K \cdot (J_n - G)$, for $G$ any positive $n$-by-$n$
matrix of rank less than or equal to $\max(\nu,n)$.

So, if $k$ has the $n$-point $\M_{1,\nu}$
Nevanlinna-Pick property, then we can choose $G$ to be the rank $(n-1)$ matrix
used in the proof of Theorem~\ref{thm}, and conclude that $F_{n+1}$ has to be
positive. But the positivity of $F_{n+1}$ clearly implies that $k$ has the
$n$-point $\M_{\mu,\nu}$ Nevanlinna-Pick property
for all values of $\mu$ and $\nu$.
\ep
\begin{cor}
The kernel $k$ is a complete Nevanlinna-Pick kernel if and only if it is a
vector-valued Nevanlinna-Pick kernel.
\end{cor}

See \cite{agmc_loc} for another approach to describing
$\M_{\nu,\nu}$ Nevanlinna-Pick kernels when there is a distinguished operator
(or tuple of operators) acting on $\H$ for which all the $k_x$'s are
eigenvectors.

\section{Representation of Complete Nevanlinna-Pick kernels}
\label{sec2}

It is a consequence of Theorem~\ref{thm} that all complete Nevanlinna-Pick
kernels have a very specific form.

\bt
\label{thm:re}
The irreducible kernel $k$ on $X$ is a complete Nevanlinna-Pick kernel if and
only if there is a positive semi-definite function $F: X \times X \to \D$
and a nowhere vanishing function $\delta$ on $X$ so that
\be
\label{rep5}
k(x,y) \= \frac{\overline{\delta(x)} \delta(y)}{1 - F(x,y)}.
\ee 
\et

\bp
(Sufficiency): If $k$ has the form of (\ref{rep5}), then $1/k$ is a rank-one
operator minus a positive operator, so has exactly one positive eigenvalue,
and the result follows from Corollary~\ref{cor:3}.

(Necessity): Suppose $k$ is a complete Nevanlinna-Pick kernel. Fix
any point $x_0$ in $X$. Then the kernel
\be
\label{eq:f}
F(x,y) \= 1 - \frac{k(x,x_0) k(x_0,y)}{k(x,y)k(x_0,x_0)}
\ee
is positive semi-definite by Theorem~\ref{thm}.
Let 
$$\delta(x) =  \frac{k(x_0,x)}{\sqrt{k(x_0,x_0)}}.
$$
It is immediate that equation~(\ref{rep5}) is satisfied.
As $k(x,x)$ is positive and finite for all $x$,  
$F(x,x)$ must always lie in $[0,1)$; as $F(x,y)$ is a 
positive semi-definite kernel, it follows that 
$|F(x,y)| < 1$ for all $x,y$.
\ep

Any positive definite kernel $k(x,y)$ can be rescaled
by multiplying by a nowhere-vanishing rank-one kernel
$\overline{\delta(x)} \delta(y)$. Let $j(x,y) =
\overline{\delta(x)} \delta(y) k(x,y)$. Then the Hilbert space
$\H_j$ is just a rescaled copy of $\H_k$: a function $f$ is in
$\H_k$ if and only if $\delta f$ is in $\H_j$, so $\H_j = \delta \H_k$.
 The multipliers of
$\H_k$ and $\H_j$ are the same, and one space has the complete
Nevanlinna-Pick property if and only if the other one does
(the matrices $F_n$ are identical, as the scaling factors cancel).
We shall say that the kernel $k$ is {\it normalized at $x_0$} if
$k(x_0,x) = 1$ for all $x$; this is equivalent to scaling the kernel by
$\displaystyle
\frac{\sqrt{k(x_0,x_0)}}{k(x_0,x)}$,
and means that in (\ref{unik}) $\delta$ can be chosen to be one,
and $F(x,y)$ becomes $\displaystyle 1 - \frac{1}{k(x,y)}$.

\section{The Universal Complete Nevanlinna-Pick Kernels}
\label{sec3}

It follows from Theorem~\ref{thm:re} that there is a {\it universal}
complete Nevanlinna-Pick kernel (actually a family of them, indexed by the
cardinal numbers).  Let $l^2_m$ be $m$-dimensional
Hilbert space, where $m$ is any cardinal bigger than or equal to $1$.
Let $\B_m$ be the unit ball in $l^2_m$, and define a kernel $a_m$ on $\B_m$ by
\be
\label{unik}
a_m(x,y) \= \frac{1}{1-\la x,y \ra}
\ee
Let $H^2_m$ be the completion of the linear span of the functions
$\{ a_m(\cdot,y) : y \in \B_m\}$, with inner product defined by
$\la  a_m(\cdot,y) ,  a_m(\cdot,x)\ra = a_m(x,y)$.
We shall show that the
 spaces $H^2_m$ are universal complete Nevanlinna-Pick spaces.

\bt
\label{thm:univ}
Let $k$ be an irreducible kernel on $X$. Let $m$ 
be the rank of the Hermitian form $F$ defined by 
(\ref{eq:f}). Then $k$ is a complete Nevanlinna-Pick kernel if and only if
there is an injective
 function $f: X \to \B_m$ and a nowhere vanishing function $\delta$
on $X$ such that 
\be
\label{eq:des}
k(x,y) \= \overline{\delta(x)}\delta(y)\, a_m(f(x),f(y)).
\ee
Moreover if this happens, then the map $k_x \mapsto 
\overline{\delta(x)} (a_m)_{f(x)}$ extends to an isometric linear embedding 
of $\H_k$ into
$\delta H^2_m$.

If in addition there is a topology on $X$ so that $k$ is continuous on $X
\times X$,
then the map $f$ will be a continuous embedding
of $X$ into $\B_m$.
\et
\bp
(Sufficiency): Any kernel of the form (\ref{eq:des}) is of the form
(\ref{rep5}).

(Necessity): Suppose $k$ is a complete Nevanlinna-Pick kernel. As $F$ is
positive semi-definite, there exists a Hilbert space of dimension $m$ (which
we shall take to be $l^2_m$) and
a map $f: X \to l^2_m$ so that $F(x,y) = \la f(x), f(y)\ra$.
Moreover, as $F$ takes value in $\D$, $f$ actuallly maps into $\B_m$.
It now follows from Theorem~\ref{thm:re} that $k$ has the form
(\ref{eq:des}).

The linear map that sends $k_x$ to the function
$\displaystyle \frac{\overline{\delta(x)}}{1-\la f(x),\cdot \ra}$
is an isometry on $\vee\{k_x : x \in X\}$ by (\ref{eq:des})
and gives the desired embedding.

If $f(x) = f(y)$ then $k_x = k_y$; as $k$ is positive definite, this
implies $x=y$.

Finally, $f$ can be realised as the composition of 
the four maps 
\begin{eqnarray*}
x &\mapsto& k_x \\
k_x &\mapsto& \overline{\delta(x)}a_m(f(x),\cdot) \\
\overline{\delta(y)}a_m(y,\cdot) &\mapsto& a_m(y,\cdot) \\
a_m(y,\cdot) &\mapsto& y
\end{eqnarray*}
The fourth map is continuous by direct calculation, the second is 
an isometry by the theorem, and the first and third maps are continuous if
$k$ is continuous.
\ep

Note that if one first normalizes $k$ at some point, $\delta$ can be 
taken to be $1$ in Theorem~\ref{thm:univ}.

For  $m=1$, the space $H^2_m$ is the regular Hardy space on the unit disk.
For larger $m$, it is a Hilbert space of analytic functions on the ball
$\B_m$. Thus every reproducing kernel Hilbert space with the complete
Nevanlinna-Pick property is a restriction of a space of analytic functions.

It was shown in \cite{ag90b} that the Sobolov space 
$W^{1,2}[0,1]$, the functions $g$ on the 
unit interval for which $\int_0^1 |g|^2 + |g'|^2 dx$
is finite, has the Nevanlinna-Pick property.
It follows from \cite[Cor. 6.5]{qui93} that the condition of 
Corollary~\ref{cor:3} is satisfied, so the Sobolov space has the complete
Nevanlinna-Pick property.
We can normalize $W^{1,2}[0,1]$ at $1$ say, by calculating that
$\displaystyle 
k_1(t) = 
{\rm cosinh}(1)\cosh(t)$,
and hence $\displaystyle
\delta(t) =
\sqrt{\sinh(1)\cosh(1)}{\rm cosech}(t)$.
Therefore
there is a continuous embedding $f:[0,1] \to \B_{\aleph_0}$ 
so that if $g$ is any function in $W^{1,2}[0,1]$, 
then $(\delta .g) \circ f^{-1}$
 extends off the curve 
$f([0,1])$ to be analytic on all of $\B_{\aleph_0}$ - even though 
$\delta .g$ need not be analytic in any neighborhood of the unit interval on
which it is originally defined.

After normalization, every separable 
reproducing kernel Hilbert space with the complete
Nevanlinna-Pick property is the restriction of the space
$H^2_{\aleph_0}$ to a subspace spanned by a set of kernel functions,
which is why we call this space universal.
The kernel $k$ is just the restriction of $a_{\aleph_0}$ to a subset of
$\B_{\aleph_0}$.

Let $\cal A$ be a normed algebra of functions on a set $X$
with the complete Nevanlinna-Pick property, \ie
there exists a positive definite function $k$ on $X \times X$ 
such that there is a function $f$ in ${\cal A} \otimes \M_k$
of norm at most one and with $f(x_i) = \L_i$ if and only if 
the  $nk$-by-$nk$ matrix
$$
k(x_i,x_j) \otimes [I_k - \Lambda_i^\ast \L_j ]
$$
is positive. It is then immediate that $\cal A$ is the multiplier algebra of
$\H_k$, and $k$ is a complete Nevanlinna-Pick kernel. If $\H_k$ is separable,
$k$ is therefore 
the restriction of $a_{\aleph_0}$ to some subset of 
$\B_{\aleph_0}$. By the Pick property, every function in $\cal A$ extends to
an element of the multiplier algebra of $H^2_{\aleph_0}$ {\it without increasing
the norm.} So every separably acting algebra with the  complete
Nevanlinna-Pick property embeds isometrically  in the multiplier algebra of
$H^2_{\aleph_0}$.

It is probably the universality of the kernel $a_m$ which is responsible for
the recent surge of interest in it - see \eg 
\cite{agmc_loc,arp97,arv97,dapi97}.

%

\end{document}